%% file: InvDensity.tex
\documentclass{article}
\usepackage{latexsym}
\usepackage{amssymb}
\usepackage{amsmath}
\usepackage{enumerate}
\usepackage[dvips]{graphicx}
\usepackage[latin1]{inputenc}
\title{Nice sets and invariant densities in complex dynamics}
\author{Neil Dobbs}

\newcommand\beginpf{\noindent \emph{Proof: }}
\newcommand\eprf{\hfill$\Box$}

\newcommand\arr{\mathbb{R}}
\newcommand\ccc{\mathbb{C}}
\newcommand\dist{\mathrm{dist}}

\newcommand\scrM{{\cal M}}
\newcommand\scrP{{\cal P}}

\newcommand\scrV{{\cal V}}

\newcommand\J{{\cal{J}}}
\newcommand\Orb{\mathrm{Orb}}
\newcommand\Swan{{{\'S}wi{\c{a}}tek}}

\newcommand\Crit{\mathrm{Crit}}

\newcommand\cbar{{\overline{\mathbb{C}}}}

\newcommand\remark{\noindent \emph{Remark: }}

\newtheorem{thm}{Theorem}
\newtheorem{dfn}[thm]{Definition}
\newtheorem{lem}[thm]{Lemma}
\newtheorem{prop}[thm]{Proposition}

\newtheorem{Fact}[thm]{Fact}

\begin{document}
\date{\today}
\maketitle
\begin{abstract}
In complex dynamics, we construct a so-called nice set (one for which the first return map is Markov) around any point which is in the Julia set but not in  the post-singular set, adapting a construction of Rivera-Letelier. This simplifies the study of absolutely continuous invariant measures. We prove a converse to a recent theorem of Kotus and \'Swi\c atek, so for a certain class of meromorphic maps the absolutely continuous invariant measure is finite if and only if an integrability condition is satisfied.
\end{abstract}
\input{body}

\bibliography{references}
\bibliographystyle{plain}

\end{document}

%% file: body.tex
\section{Introduction}
In dynamical systems, invariant measures which are absolutely continuous with respect to some natural reference measure (often Lebesgue measure, but also other \emph{conformal} measures) are of great interest. We shall study such measures in the setting of complex one-dimensional dynamics. Both rational and transcendental maps are considered. In the transcendental setting one must deal with maps with unbounded derivative.

This work starts with a useful construction (in Proposition \ref{prop:nice}) of simply connected \emph{nice sets} around each point in the Julia set but not in the post-singular set. First return maps to such maps have good Markov properties. We use this to give a very simple proof of a couple of essentially known results concerning existence and smoothness of invariant densities (Theorem \ref{thm:BockConf}). Then we give a converse to a theorem of Kotus and \'Swi\c atek (\cite{KotSwi:Mero}): for a class of meromorphic maps admitting absolutely continuous invariant measures (\emph{acims}), the acims are finite if and \emph{only if} an integrability condition is satisfied (Theorem \ref{thm:converse}).  

We shall introduce some notation and definitions now, before discussing the results in more detail in subsections \ref{sec:niceresults}-\ref{sec:converse}.

Let $f : \ccc \to \cbar$ be a transcendental meromorphic function or let $f: \cbar \to \cbar$ be a rational map of degree $\geq 2$. The Fatou set is defined as usual using normal families: a point $z$ is in the Fatou set if and  only if there is a neighbourhood of $z$ on which the iterates of $f$ are well-defined and form a normal family. The Julia set, $\J(f)$, of $f$ is the complement in $\cbar$ of the Fatou set. We use the spherical metric. 
The $\omega$-limit set of a point $z$ is denoted $\omega(z) \subset \cbar$. Given a point $z \in \cbar$, we set 
$$z \in \Orb(z):= \{f^i(z) : i \geq 0 \mbox{ and }f^i(z)\mbox{ is defined} \}.$$ Similarly, for a set $A$, we write $\Orb(A) := \bigcup_{z \in A} \Orb(z)$.  

Let $S(f) \subset \cbar$ denote the set of \emph{singular values} of $f$: a point $z$ is an element of $\cbar \setminus S(f)$ \emph{if and only if} there is a neighbourhood of $z$ on which all inverse branches are well-defined, univalent maps. $S(f)$ contains all critical and asymptotic values. 
Set 
$$
\scrP(f) := \overline{\Orb(S(f))}. 
$$
$P(f)$ is called the \emph{post-singular set}.

All measures considered shall be assumed Borelian without further reference. We denote by $\scrM_\infty(f)$ the set of all conservative, $\sigma$-finite, ergodic $f$-invariant measures. We set 
$$\scrM(f):= \{ \mu \in \scrM_\infty(f) : \mu(\cbar) = 1\}.$$ 

If the Fatou set is non-empty, or if one wants to consider the pressure function and thermodynamic formalism, then conformal measures are of interest.
\begin{dfn} Let $p, t \in \arr$. 
    We call a measure \emph{$(t,p)$-conformal}, with respect to some metric $\rho$, if 
    \begin{itemize}
        \item
            it is finite on compact sets disjoint from $\scrP(f)$; 
        \item
            it has Jacobian $\exp(p)|Df|_\rho^t$,  wherever this is finite, where $|Df|_\rho$ denotes the norm of the derivative with respect to that metric, and with the convention that $0^0 =1$;
        \item
            if $|Df(x)|_\rho^t = \infty$, then $m(x) = 0$ and $m(f(x))$ is zero or finite. 
    \end{itemize}
\end{dfn}
With this definition, if $t >0$ then critical values have measure zero, critical points need not have. If $t < 0$ then critical points have measure zero, critical values need not have. 
Standard and spherical Lebesgue measures are $(2,0)$-conformal with respect to the Euclidean and spherical metrics respectively. 

If the metric is unspecified, the norm of the derivative $|Df|$ shall henceforth be taken with respect to the spherical metric. 

\subsection{Nice sets}\label{sec:niceresults}
\begin{dfn} An open set $U$ is called \emph{nice} if $f^n(\partial U) \cap U =
\emptyset$ for all $n >0$.
\end{dfn}
This implies that every pair of pullbacks of a nice set $U$
(connected components of $f^{-n}(U), f^{-n'}(U)$ for some $n,n'\geq 0$)
is either nested or disjoint.

In the Section \ref{sec:nice} we adapt a construction of nice sets due to J. Rivera-Letelier (\cite{Rivera:Connecting, PRL:StatTCE}) to prove the following result. In the setting of real interval dynamics, existence of nice sets (or \emph{regularly returning} sets) is easy to prove using (pre-)periodic points. In higher dimensions the existence of such sets is highly non-trivial. 

\begin{prop}\label{prop:nice}
Suppose $z\in \J(f)\setminus \scrP(f)$.
Then $z$ is contained in a simply connected nice set $U_z$ with arbitrarily small diameter satisfying $\dist(U_z, \scrP(f)) > 0$. 
\end{prop}

\subsection{Existence and smoothness of invariant densities}
Call a point $z$ \emph{univalently inaccessible} if there is an open set $U$ with non-empty intersection with the Julia set, such that there does \emph{not exist} a triple $(y,n,m)$ satsifying the following:
\begin{itemize}
    \item
        $y \in U$;
    \item
        $n,m \geq 0$;
    \item 
        $f^n(z) = f^m(y)$;
    \item
        $0 < |Df^n(z)||Df^m(y)| < \infty$.
\end{itemize}
Let $E$ denote the set of all univalently inaccessible points. 
For most maps $E\cap \J(f) = \emptyset$. 
The simplest map for which $E\cap \J(f)$ is non-empty is the Chebyshev polynomial $f:z\mapsto 4z(1-z)$, for which $E = \{0, 1, \infty\}$.

\begin{thm} \label{thm:BockConf}
Let $f : \ccc \to \cbar$ be a meromorphic function. 
Let $p,t \in \arr$ and suppose $m$ is a $(t,p)$-conformal measure with respect to either the Euclidean metric or the spherical metric.  
Suppose
$m\left(\{ z \in \J  : \omega(z) \not\subset \scrP(f)\}\right) > 0$ and $m(E)=0$.

Then $f$ admits a non-atomic measure $\mu \in \scrM_\infty(f)$, equivalent to $m$ with the support of $\mu$ equal to the Julia set. 
On a neighbourhood of any point in $\J \setminus \scrP(f)$, 
there is a density function $\frac{d\mu}{dm}$ which 
is analytic and non-zero. Moreover, $m(\scrP(f)) = 0$. 
\end{thm}

\remark If $\J = \cbar$, then it follows immediately that the density is analytic on $\J \setminus \scrP$. In any case, this local result is sufficient for most purposes.


Theorem \ref{thm:BockConf} will not surprise the specialist. For rational maps, when the reference measure is Lebesgue measure, ergodicity and conservativity of Lebesgue measure date back to the work of Lyubich \cite{Lyubich:Typical}. Grzegorczyk \emph{et al}, in \cite{GPS:Rational}, give a detailed proof of the existence of the absolutely continuous invariant probability measure. For meromorphic maps,   ergodicity and conservativity of Lebesgue measure were shown by H. Bock in \cite{Bock, Bock2}. Kotus and Urbanski (\cite{KotUrb03}) showed that in fact there is a $\sigma$-finite invariant measure. In 
\cite{przytycki1999rtr}, Przytycki and Urba\'nski showed something similar, including analyticity of the density, in the rational setting  for $(t,0)$-conformal measures. They also discuss what happens when the Julia set is contained in a finite union of \emph{real-analytic} sets.  
In \cite{KotSwi08}, it is shown that, if the acim with respect to Lebesgue measure is a probability measure, then it has density bounded from below on an open set. Compare \cite{MeBS08}, where a predecessor of Theorem \ref{thm:BockConf} can be found.

 However, the union of these proofs sprawls unnecessarily and we wish to give a reasonably short and elegant proof, also in Section \ref{sec:nice}, of a more general result covering both the rational and transcendental cases.

 \subsection{Finiteness of acims} \label{sec:converse}
 In this subsection, we only use the Euclidean metric. 
\begin{thm}\label{thm:converse} 
Let $f : \ccc \to \cbar$ be a meromorphic function such that 
there exists a positive Lebesgue measure set of points $z \in \J(f)$ such that $\omega(z) \not\subset \scrP(f)_*$. 

Let $A$ be a forward-invariant, bounded set and suppose
$f$ admits  
 a pole of order $M$ which is not an omitted value.  
If 
the $\sigma$-finite measure given by Theorem \ref{thm:BockConf} is finite, then 
\begin{equation}\label{eqn:est}
\int_{|z| > r_0} \frac{-\log \dist(f(z),A)}{|z|^{2+\frac{2}{M}}} dm  < \infty
\end{equation}
for some $r_0 > 0$, where integration is with respect to Euclidean Lebesgue measure. 
\end{thm}
\remark By Theorem \ref{thm:BockConf}, $\J(f) = \cbar$. One can give a version of this theorem where the reference measure is $(t,p)$-conformal rather than Lebesgue measure. One must replace $(2+2/M)$ by 
$(t+t/M)$, then, provided the conformal measure has support equal to $\cbar$, the proof is the same. We choose not to do so in order to keep the discussion concise. 

\remark Suppose $f$ admits an asymptotic value whose orbit is bounded. Then (\ref{eqn:est}) implies
$$
\int_{|z| > r_0} \frac{-\log \dist(f(z),a)}{|z|^{2+\frac{2}{M}}} dm  < \infty
$$
 One can rewrite the inequality as 
$$
\int_{r_0}^\infty \frac{m(r,a)}{r^{1+ \frac{2}{M}}} dr < \infty
$$ 
where 
$$
m(r,a) = \int_0^{2\pi} \log^+ \frac{1}{\dist(f(r e^{i\theta}), a)} dr.
$$
The notation $m(r,a)$ for this quantity is from Nevanlinna theory (see \cite{Tsuji:Book}, \cite{KotSwi:Mero}). 

Theorem \ref{thm:converse} is new. Note that it is only interesting for transcendental maps, since in the rational setting the integral is always finite. Recently Kotus and \Swan \ proved the following  nice result which motivated our theorem.
\begin{Fact}[\cite{KotSwi:Mero}]  \label{Fact:KS}
Let $m$ denote Lebesgue measure.
Let $f : \ccc \to \cbar$ be a meromorphic map with finitely many singular values. Suppose all poles of $f$ have multiplicities bounded by $M$. Suppose also that $\J(f) = \cbar$ and $\scrP(f) \cap (\Crit(f) \cup \{\infty\}) = \emptyset$ and that,  for some $r_0 >0$,
$$
\int_{r_0}^\infty \frac{m(r,a)}{r^{1+\frac{2}{M}}} dr < \infty
$$
for each asymptotic value $a$.

Then there exists $\mu \in \scrM(f)$ which is absolutely continuous with respect to Lebesgue measure.
\end{Fact}
They also provide an example of a map satisfying all of the assumptions bar the integrability one and show that this map does not admit a \emph{finite}, absolutely continuous, invariant measure. Theorem \ref{thm:converse} states that this integrability condition 
(and indeed a stronger one, integrability of the logarithmic distance to any forward invariant, bounded set) is necessary for a finite probability measure to exist for a whole class of maps including those considered in  Fact \ref{Fact:KS}. 

One could be tempted to view  Fact \ref{Fact:KS} as an analogue to the equivalent result for Misiurewicz rational maps. However, compare it with the following result of Benedicks and Misiurewicz which we would view as more relevant:
\begin{Fact}[\cite{MisBen:Flat}] \label{Fact:MisBen}
Let $f : I \to I$ be a $C^3$ map of the interval $I$ with a unique critical point $c$. Suppose that $c \notin \omega(c)$ and that all periodic points are repelling. Then there exists a $\sigma$-finite, conservative, ergodic, invariant measure $\mu$ which is absolutely continuous with respect to Lebesgue measure. 
The measure $\mu$ is finite if and only if
$$
\int \log |f'(x)| dx > -\infty.$$
\end{Fact}
They also showed that this inequality is equivalent to 
$$
\int -\log \dist(f(x),f(c)) dx < \infty,$$
strikingly similar to Fact \ref{Fact:KS}. The reason the interval setting is perhaps more similar is that rational maps are very rigid: the critical points are all like $z^n$ for some $n$, so they cannot scrunch up space too much. On the interval, one can have flat critical points of type $\exp(|x|^{-\alpha})$ for $\alpha >0$. 
Fact \ref{Fact:MisBen} implies that for the maps considered with a critical point of this type, the measure is finite if and only if $\alpha < 1$. For meromorphic maps, large regions near infinity can get mapped very close to asymptotic values. Theorem \ref{thm:converse} roughly says that if (Lebesgue-) many points get mapped  too close to an asymptotic value and then take too long to leave a bounded region, this is an obstruction to the measure being finite. 

In prior work \cite{Me:Cusp}, we generalised one direction of Fact \ref{Fact:MisBen}: 
\begin{Fact}[\cite{Me:Cusp}] \label{Fact:MeCusp} Let $f : I \to I$ be a $C^{1+\epsilon}$ map of the interval $I$. Suppose there exists $\mu \in \scrM(f)$ such that the entropy of $\mu$ is positive and that $\mu$ is absolutely continuous with respect to Lebesgue measure. Then the support of $\mu$ is a finite union $X$ of closed intervals on which
$$
\int_X \log |f'(x)| dx > -\infty.$$
\end{Fact}
We view Theorem \ref{thm:converse} not only as a converse to Fact \ref{Fact:KS},  but also as a first generalisation of Fact \ref{Fact:MeCusp}
 in the meromorphic setting.

\section{Nice sets and invariant densities} \label{sec:nice}

Our proof of Theorem \ref{thm:BockConf} relies on an adaptation of the construction of \emph{nice sets} due to J. Rivera-Letelier (\cite{Rivera:Connecting, PRL:StatTCE}). His construction is around critical points and uses different methods to guarantee decrease of preimages along backward branches. We use nice sets to find induced Markov maps on a neighbourhood of each point of $\ccc \setminus \scrP(f)$. Theorem \ref{thm:BockConf} will then follow by standard arguments. 

Recall that we use the spherical metric, not the Euclidean one. 

For $n \geq 0, z\in \cbar, 0<R$ and $U$ an open set compactly contained in $B(z,R)$, let $\scrV_n(U,z,R)$ denote the collection of open sets $V$ such that $f^n (V)= U$, $V \subset V'$ for some 
conected $V'$ on which $f^n$ is univalent and $f^n(V') = B(z,R)$. Let $\scrV(U,z,R) := \bigcup_{n>0} \scrV_n(U,z,R)$.

The following lemma is standard, and very similar to Lemma 1 of \cite{GPS:Rational}, for example.
\begin{lem} \label{lem:decrease} Let $z \in \J(f), R > 0$ and $U$ be an open set compactly contained in $B(z,R)$. 
For each $K>0$ there exists an $N \geq 1$ such that if $n \geq N$ and $V \in \scrV_n(U,z,R)$, then $|Df^n(z')| > K$ for all $z' \in V$.
\end{lem}
\beginpf
There are infinitely many periodic orbits (we assume $f$ is not a Möbius transformation (in particular $\J\ne \emptyset$)), so we can fix three points whose orbits are disjoint from a small neighbourhood $W$ of $z$ with $W \subset B(z,R)$. 
Then the holomorphic inverse branches of $f^{-n}_{|W}$ omit the three points, so they form a normal family by Montel's Theorem. 

By the Koebe principle, there is a uniform bound on the distortion of $f^n$ on each $V \in \scrV_n(U,z,R)$ for all $n$.
Suppose such an $N$ (from the statement) does not exist for some $K>0$. Then, by uniformly bounded distortion, there is a sequence of holomorphic inverse branches $g_n$ on $B(z,R)$ with derivative at $z$ uniformly bounded away from zero. 

Then, by Montel,  there is a subsequence which converges to a non-constant function $g$ on a neighbourhood of $z$. 
Then a neighbourhood of $g(z)$ is mapped infinitely often by iterates of $f$ onto a neighbourhood of $z$ and into $W$. It follows that $z \notin \J(f)$, contradiction. 
%
\eprf

Now we need to show that nothing too nasty happens for small $n$. The complication is that $z$ could be in the backward orbit of a pole, and so be an accumulation point for the set $f^{-k}(z)$ for some fixed $k$. If this were impossible, the following lemma would be simpler. 
\begin{lem} \label{lem:branches}
    Let $z \in \cbar$, $R >0$, $\tau >0$ and $N \geq 1$. Let $r < R$. Then there exists $M>0$ so that, for $n \leq N$ all holomorphic inverse branches $g$ of $f^{-n}$ on $B(z,R)$ satisfy $|Dg| < M$ on $B(z,r)$. Moreover, all but a finite number of such inverse branches satisfy $|Dg| < \tau$ on $B(z,r)$.
\end{lem}
\beginpf
The key is that the images of $B(z,r)$ by any pair of inverse branches of $f^{-n}$ are pairwise disjoint. Thus there are at most a finite number of any minimal size (Koebe implies the images are not too distorted, and the Riemann sphere has finite area), so firstly there is an upper bound to the size (and, equivalently, derivative), and secondly all but a finite number of such inverse branches have small images. If the image is small, the derivative $|Dg|$ is small. 
\eprf


\begin{lem} \label{lem:knice} Let $z \in \J(f)$ and $R >0$. Suppose $z$ is not a parabolic periodic point. Let $\kappa > 1$. For all $r>0$ sufficiently small, there exists an open, simply connected set $U = U_r(z)$  such that
\begin{itemize}
\item
$B(z,r) \subset U \subset B(z,\kappa r)$;
\item
if $n > 0$, 
$V \in \scrV_n(U,z,R)$  and $V \cap U \ne \emptyset$ then $V \subset U$;
\item there exists $\theta >1$ so that, for all $n$, if  $V \in \scrV_n(U,z,R)$ and $V\cap U \ne \emptyset$, then $|Df^n(z')| > \theta$ for all $z' \in V$.
\end{itemize}
\end{lem}
\beginpf 
%
We call an inverse branch $g$ \emph{$z$-periodic} if $g(z) = z$. This of course can only happen if $z$ is periodic. 

By Lemma \ref{lem:decrease}, there exist $\delta$ and  $N$ such that for all $n \geq N$ and all $W \in \scrV_n(B(z,\delta),z,R)$,  $|Df^n(z)| > 2\kappa/(\kappa-1)$ for all $z \in W$. 
Given a holomorphic non-$z$-periodic inverse branch $g$ of $f^n$ on $B(z,R)$, if $r$ is sufficiently small then $g(B(z,r)) \cap B(z,\kappa r) = \emptyset$. This allows one to disregard a finite number of branches to deduce the following, by 
Lemma \ref{lem:branches}: There exists $0< r < \delta$ sufficiently small that, if $n < N$,  any non-$z$-periodic holomorphic inverse branch $g$ of $f^n$ on $B(z,R)$  such that $g(B(z,r)) \cap B(z,\kappa  r) \ne \emptyset$ satisfies 
$$
|Dg| < (\kappa-1)/2\kappa
$$
 on $B(z,r)$. 
In particular, if $W \in \scrV_n(B(z,\kappa r),z,R)$ for some $n > 0$ and $W \cap B(z,\kappa r) \ne \emptyset$ then either $|Df^n| > 2\kappa/(1-\kappa)$ on $W$ or $z$ is periodic and $z \in W$. 

If $z$ is periodic we also require $r$ to be small enough that $|Df| > 1$ on $B(z,\kappa r)$.

Set $U_0 := B(z,r)$. 
For $n > 0$ define $U_n$ as the connected component of
$$
U_0 \cup \bigcup_{i=1}^n \bigcup_{W \in \scrV_i(B(z,r),z,R)} W$$
containing $z$ and $U_r'(z) := \bigcup_{n \geq 0} U_n$. We
prove by induction that $U_n \subset B(z,\kappa r)$ for all $n \geq
0$. This is clearly true for $n = 0$. So suppose it is true for all $n
\leq k$. We must show it holds for $n = k+1$.

Let $X$ be a connected component of $U_{k+1} \setminus U_0$. Then
there is a minimal $m \geq 0$ such that there is a $y\in X$ and $V \in \scrV_m(B(z,r),z,R)$. Let $g$ be the corresponding inverse branch, so $f^m(y) \in U_0$ and $g(f^m(y)) = y$.
 Since $y \notin U_0$, $m \geq 1$ and  $g$ is not $z$-periodic.  Consider $f^m(X)$. This set is
contained in $U_{k+1-m}$, and so by hypothesis is contained in
$B(z,\kappa r)$. But then $X$ is contained in $g(B(z,\kappa r))$
 so $|X| < ((\kappa-1)/2\kappa) |B(z,\kappa r)| = r(\kappa-1)$. Thus $X \subset B(z,\kappa r)$ as required, completing the inductive argument. 

 Fill in the set $U'_r(z)$ to get a simply connected set $U = U_r(z)$ (the smallest simply connected set containing $U_r'(z)$ and contained in $B(z, \kappa r)$).

Now suppose there exists a $V \in \scrV_n(U,z,R)$ and a point $y \in  V \cap \partial U$. Then $f^n(y)$ is in some $W \in \scrV_m(B(z,r),z,R)$. But this means $y$ is contained in some element of $\scrV_{n+m}(B(z,r),z,R)$ which must intersect $U$ since $y\in \partial U$, so $y$ is contained in $U$, contradiction. 
\eprf 
 
Note that one can construct a similarly defined set around parabolic points if one excludes periodic inverse branches too. 

\bigskip

\noindent \textbf{Proof of Proposition \ref{prop:nice}:}
Recall that an open set $U$ is called \emph{nice} if $f^n(\partial U) \cap U =
\emptyset$ for all $n >0$. 
Suppose $z\in \J(f)\setminus \scrP(f)$. Then $z$ is not parabolic. Take some small $R < \dist(z,\scrP(f))$. 
Then all inverse branches on $B(z,R)$ are holomorphic.
Let $U$ be the set given by 
  Lemma \ref{lem:knice}.  We claim this set is nice. Indeed, suppose $x \in \partial U$ and $f^n(x) \in U$. Then $x \in V$, for some $V \in \scrV_n(U, z, R)$. But $V \subset U$ so $x \in U$, contradiction. 
\eprf

Given a nice set $U$, we write $r_U(z) :=\inf\{k \geq 1: f^{r_U(z)}(z) \in U\}$ for the first return time to $U$ and define the first return map $\phi_U$ by 
$$
\phi_U(z) := f^{r_U(z)}(z)
$$
for $z \in U$ such that $r_U(z)$ is defined (finite). By the nested or disjoint property, the domain of definition of $\phi_U$ is a countable union of connected components $U_i$ on which $r_U$ is constant. If $U$ is connected then 
$$
\phi_U : U_i \to U
$$
is biholomorphic for each $i$. Let $m$ be a $(t,p)$-conformal measure for $f$. Then  the Jacobian of $\phi_U$ at $z$ is $\exp(p r_U(z)) |D\phi_U(z)|_\rho^t$.

\begin{lem}[Folklore Theorem]\label{Fact:Folklore}
Let $U$ be a bounded, simply-connected nice set such that $\dist(U, \scrP(f)) > 0$. Let $\Lambda_U :=\bigcap_{k\geq0} \phi_U^{-k}(U)$. If $m(\Lambda_U) >0$ then there exists a non-atomic, ergodic, absolutely continuous, $\phi_U$-invariant, probability measure $\nu$ with 
a density $\rho_\nu$ which is non-zero and analytic on $U$. The measure $\nu$ is the unique absolutely continuous invariant probability measure on $U$. 
Moreover, $m(U \setminus \Lambda_U)=0$ and $m$-almost every point in $U$ is recurrent. 
\end{lem}
\beginpf 
%
%
This is standard, see for example Theorem 6.1.3 in \cite{MauldinUrbanski:Book}.
That the density is smooth is easy, the argument for analyticity can also be found in \cite{przytycki1999rtr} or \cite{MPU:Rigidity}.
To deduce $m(U \setminus \Lambda_U)=0$, note, by bounded distortion and conformality of the measure, there exists $C>0$ such that, for any two subsets $A,B\subset U$ and $k\geq0$, 
$$
C^{-1}\frac{m(A)}{m(B)} \leq \frac{m(\phi_U^{-k}(A))}{m(\phi_U^{-k}(B))} \leq C\frac{m(A)}{m(B)}.
$$
Let $A_n$ denote the set of points $x\in U$ for which $\phi^n_U(x)$ is defined, but $\phi^{n+1}_U(x)$ is not. Then $U = \Lambda_U \cup \bigsqcup_{n\geq0}A_n$, and $\phi_U^{-1}(A_n) = A_{n+1}$. Setting $A := A_0$ and $B:=\Lambda_U$ in the above inequality, we have $A_n \geq cA_0$ for some constant $c>0$ and all $n$. Since the measure is finite on $U$, it is zero on each $A_n$, as required. 

That $\nu$ is non-atomic follows by a similar argument. Remark just that $\phi$ has more than one branch (an infinity of branches, actually), so that if $m$ admits an atom, then it has a non-periodic atom at a point $q$; set $A_0 := \{q\}$ and $A_n := \phi_U^{-n}(A_0)$. 

Recurrence follows from ergodicity.
\eprf

%
\smallskip

\noindent \textbf{Proof of Theorem \ref{thm:BockConf}.}  
Consider  the compact sets 
$$C_k :=\ccc \setminus\left(B\left(\scrP(f),\frac{1}{k}\right) \cup \{z:|z| >k\}\right).
$$

Recall that, according to the hypotheses, there is a positive measure set of points $z$ with $\omega(z) \not\subset \scrP(f)$. But then there is some $C_k$ and a positive measure set of points $z$ with $\omega(z) \cap C_k \ne \emptyset$. Each point in $C_k$ is contained in a nice set given by Proposition~\ref{prop:nice}, so there is a finite cover of $C_k$ by such sets, and thus there is one, $U$ say, and a set $A$ of positive measure such that, for all $z \in A$, $\omega(z) \cap U \ne \emptyset$. 

We can write $A = \bigcup_{n \geq 0} A_n$ where $A_n = \{y \in A : r_{U}(y) = n\}$. Then there is some  $A_n$ 
of positive measure. 
Then  $f^n(A_n)$ has positive measure and  
$f^n(A_n) \subset \bigcap_{k\geq 0} \phi_{U}^{-k}(U)$. We can apply the Folklore Theorem. 

By the Folklore Theorem, $m$-almost every point in $U$ is recurrent. 
We have $\scrP(f) \cap U = \emptyset$ and $\scrP(f)$ is forward-invariant.
Points in $\scrP(f)\setminus E$, are, by definition, not univalently inaccessible. By conformality, it follows then that $\scrP(f)\setminus E$ has zero measure. By hypothesis, $m(E)=0$, so $m(\scrP(f)) = 0$. 


Let $\nu$ be the $\phi_U$-invariant measure given by the Folklore Theorem. 
Let $\{U_j\}_j$ denote the connected components of the domain of $\phi_U$, and $r_j$ the first return time to $U$ on $U_j$. Now   
$$
\mu := \sum_j \sum_{k=0}^{r_j-1} f^k_*\nu_{|U_j}
$$
is in $\scrM_\infty(f)$ and absolutely continuous and equivalent to $m$.

Let $y \in \J \setminus \scrP(f)$. We consider a nice set $U_y \ni y$ given by Proposition~\ref{prop:nice} and the first return map $\phi_{U_y}$. Since $\mu/\mu(U_y)$ is an ergodic, absolutely continuous, invariant, probability measure for $\phi_{U_y}$, its density is analytic on $U_y$ (using the Folklore Theorem). Thus the density of $\mu$ is analytic on a neighbourhood of $y$. 
\eprf

\section{Finite mass}
Now let us prove Theorem \ref{thm:converse}. For these final paragraphs we use the Euclidean metric. 

Let $p$ be a pole which is not an omitted value. There exists a point $y \in \ccc$ and an $n \geq 0$ such that $f^n(y) = p$ and $y \notin \scrP(f)$.  
By Theorem \ref{thm:BockConf}, the density $\rho$ of the absolutely continuous invariant measure $\mu$ 
is bounded away from zero on a neighbourhood of $y$. Thus $\rho$ is bounded away from zero on a neighbourhood of $p$. If $p$ has multiplicity $M > 0$, then it follows that there exists $c,r_0 > 0$ such that 
$$
\rho(z) > c \frac{1}{|z|^{2+\frac{2}{M}}}
$$
for all $z$ satisfying $|z| \geq r_0$. 

Let $A$ be any subset of $\ccc$ such that $\Orb(A)$ is bounded. Then there exists $\varepsilon >0$ such that $B(\Orb(A),2\varepsilon)$ contains no poles. 
The derivative $|f'|$ is bounded from above by a positive constant $K > 1$ on $B(\Orb(A),\varepsilon)$. Now let $x \in \ccc$. 
Denote by $n(x)$ the least integer $n \geq 1$ such that $\dist(f^n(x), \Orb(A))  > \varepsilon$. Then $\varepsilon <  K^{n(x)} \dist(f(x),\Orb(A))$, so 
$$
n(x) > (1/\log K) (\log \varepsilon - \log \dist(f(x), \Orb(A))),$$
so $n(x) > - c_1 \log \dist(f(x),A) -c_2$ for some constants $c_1, c_2 > 0$.

Consider $r_0$ as before, but large enough that $\Orb(A) \subset B(0, r_0 - \varepsilon)$. 
We shall consider the first return time $r_W$ to $W := \{|z| \geq r_0\}$. 
By Kac's Lemma, $\int_W r_W(z) d\mu(z) = \mu(\ccc)$. Combined with the above density estimate,  
$$ \int_W r_W(z) \frac{1}{|z|^{2+\frac{2}{M}}} dm < \mu(\ccc)/c.$$
For $z \in W$, 
$$
r_W(z) \geq  n(f(z)) > 
- c_1 \log \dist(f(z),A) -c_2$$
so, if $\mu$ is finite, then 
$$ 
+\infty > \int_W \frac{-\log \dist (f(z),A) }{|z|^{2+\frac{2}{M}}} dm $$
as required.
\eprf

\section*{Acknowledgments}
The author is grateful to IMPAN in Warsaw where this work was started, supported by the  EU training network ``Conformal Structures and Dynamics''. The author  currently enjoys the support of the 
Göran Gustafssons Stiftelse and the Knut och Alice Wallenbergs Stiftelse. Thanks are due to a first referee and to Bartek Skorulski for helpful conversations and remarks.